\documentclass[a4paper,11pt]{article} 

\usepackage{amssymb,amsmath,amsthm}
\usepackage[T1]{fontenc} 

\usepackage[mathcal]{euscript}      
\usepackage{amscd}
\DeclareMathAlphabet{\mathpzc}{OT1}{pzc}{m}{it}

\usepackage{geometry,float}      
\usepackage[english,french]{babel}

\newtheorem{theorem}{Theorem}[section]
\newtheorem{lemma}[theorem]{Lemma}

\newtheorem{definition}[theorem]{Definition\rm}
\newtheorem{remark}{\it Remark\/}

\title{Majoration of the dimension of the space of  concatenated solutions of a specific pantograph equation}
\author{Jean-Francois Bertazzon\footnote{ Lyc\'ee Notre Dame de Sion, Marseille, France.
    Email: \texttt{jeffbertazzon@gmail.com}}}

\begin{document}

\maketitle
\selectlanguage{english}

\begin{abstract}
For each $\lambda \in \mathbb N^*$, we consider the integral equation:
\[
\int_{\lambda y} ^{\lambda x} f(t)\, d t = f(x) - f(y)  \mbox{ for every $(x,y)\in {\mathbb R}_+^2$,}
\]
where $f$ is the concatenation of two continuous functions $f_a,f_b:[0,\lambda] \rightarrow {\mathbb R}$
along a word $u= u_0u_1\cdots\in\{a,b\}^{\mathbb N}$ such that $u=\sigma(u)$,
where $\sigma$ is a $\lambda$-uniform substitution satisfying  some combinatorial conditions.

There exists some non-trivial solutions (\cite{BD}).
We show in this work that the dimension of the set of solutions  is at most two.
\end{abstract}

\section{Introduction}

For each positive integer $\lambda\geq 2$ and each integer $\delta \in {\mathbb Z}^*$, we consider the integral equation:
\begin{equation} \int_{\lambda y} ^{\lambda x} f(t)\, d t = \delta\big( f(x) - f(y)\big)  \mbox{ for every $(x,y)\in ({\mathbb R}_+)^2$.} \tag{$E_{\lambda,\delta}$} 
\label{eq:E_lambda_delta} \end{equation}
This equation is a particular case of the \emph{pantograph} equation:
\[
f'(x) = a f(\tau x) + b f(x) \quad \text{with $(a,b) \in \mathbb R^2$ and $\tau \in \mathbb R_+$ for $x\geq 0$.}
\]
We refer to \cite{panto3}, \cite{panto4}, \cite{panto1} and \cite{panto2} for more details on the pantograph equation.

We prove in \cite{BD} that we can extend each continuous function $f$ defined on $[1,\lambda]$ such that $f^{(n)}(1)=f^{(n)}(\lambda)=0$ for every non-negative integer $n$, into a continuous solution of \eqref{eq:E_lambda_delta}.
Therefore the set of continuous solutions of \eqref{eq:E_lambda_delta} is an infinite-dimensional vector space.

Moreover, we prove in \cite{BD} that the non-identically zero solutions are not periodic.
It seems natural to look for the simplest solutions of \eqref{eq:E_lambda_delta}.
The periodic functions are the repetition of the same motif.
We study the functions which are the repetition (not periodically) of two functions.
This leads us to the following notion of \emph{concatenation} of two functions along a word. 

\begin{definition} \label{conca}
Let $\lambda\geq 2$ be a positive integer and $f_a,f_b:[0,\lambda]\rightarrow {\mathbb R}$ be two functions.
For each finite word $u=u_0\cdots u_{n-1}\in \{a,b\}^n$ of length $n$, we define a function $f_u:[0,n\lambda]\rightarrow \mathbb R$  
called \emph{the concatenation of} $f_a$ and $f_b$ along $u$ by:
\[ f_u(x+\lambda k) := f_{u_k}(x)  \mbox{ for $x\in[0,\lambda]$ and $k\in \{0,\ldots,n-1\}$.} \] 
We extend this definition to  infinite words.
\end{definition} 
{}

Our main result is the following theorem. We recall in Section \ref{combi} some notions of combinatorics on words requisite to fully understand this result.
{}

\begin{theorem} \label{theorem}
We consider  a $\lambda$-uniform substitution $\sigma$, satisfying some combinatorial conditions (Relations \eqref{eqlambda} and \eqref{eqdelta}) and  $u= u_0u_1\cdots\in\{a,b\}^{\mathbb N}$ an infinite word such that $u=\sigma(u)$.

We consider the integral equation:
\begin{equation} \int_{\lambda y} ^{\lambda x} f(t)\, d t = f(x) - f(y)  \mbox{ for every $(x,y)\in ({\mathbb R}_+)^2$.} \tag{$E_{\lambda}$} \label{eq:E_lambdq} \end{equation}

We denote by $\mathcal S_\lambda$ the set of solutions $f$ of  \eqref{eq:E_lambdq} which are the concatenation of two continuous functions $f_a,f_b:[0,\lambda] \rightarrow {\mathbb R}$ along the word $u$.
Then $\mathcal S_\lambda$ is a vector space of dimension at most $2$.
\end{theorem}
{}

We prove in \cite{BD} that $\mathcal S_\lambda$ is of dimension at least $1$.
To construct a non-trivial solution, we renormalized some iterated Birkhoff sums.
The technique used to prove Theorem \ref{theorem} (in Section \ref{proof}) is very different.
It is based on the relation between the values taken by the functions and their moments.
This brings us back to the historical first non-trivial solution associated to the Prouhet-Thue-Morse substitution ($a\to ab$ and $b\to ba$) constructed by Fabius (\cite{Fabius66}) as a cumulative distribution function. 
 
We do not have examples of substitutions for which the dimension of $\mathcal S_\lambda$ is two.
\medskip

We will use the two following basic results (see \cite{BD}).
{}

\begin{remark} \label{rem1}
Let $f$ be as in Theorem \ref{theorem}, then for every finite word $v$ of length $n$:
\[ \int_{\lambda y} ^{\lambda x} f_{\sigma(v)} (t)\, d t = f_v(x) - f_v(y)  \mbox{ for every $(x,y)\in [0,n\lambda]^2$.} \]
\end{remark}
{}

\begin{remark} \label{rem2}
We have $f_a(0)=f_a(\lambda)=f_b(0)=f_b(\lambda)$.
\end{remark}

\section{Some notions about combinatorics on words} \label{combi}

We consider the \emph{alphabet} $\{a,b\}$ consisting of two \emph{letters} $a$ and $b$. 
We denote by $\{a,b\}^*$ the set of \emph{finite words}. 
Endowed with the concatenation, it is a free monoid and an endomorphism is called a \emph{substitution}.
If $u$ is a finite word, we  denote by $|u|$ its \emph{length} and $|u|_\alpha$ the number of occurrences of the letter $\alpha$ for $\alpha \in \{a,b\}$.

A substitution $\sigma$ is said to be \emph{$\lambda$-uniform} if $\lambda := |\sigma(a)|=|\sigma(b)|$.
We only consider  $\lambda$-uniform substitutions $\sigma$ such that:
\begin{equation} \label{eqlambda}
\lambda_a := |\sigma(a)|_a = |\sigma(b)|_a \quad \mbox{and} \quad \lambda_b := |\sigma(a)|_b = |\sigma(b)|_b.
\end{equation}
We have  of course $\lambda_a+\lambda_b=\lambda$.
The next notion  takes care of the order of apparitions of the letters in $\sigma(a)$ and $\sigma(b)$.
If $u=u_0\cdots u_{n-1}$ is a finite word with $n>1$, we define the  set of strict prefixes by  $\mbox{pref}(u):=\{u_0\cdots u_k ;  0\leq k < n-1\}$ and
\begin{equation} \label{eqdeltaab} 
\delta_{a  \in  \sigma(a) } ^{(1)} :=  \sum_{v \in \mbox{\small pref} (\sigma(a))} |v|_a  \mbox{ and }  \delta_{a  \in  \sigma(b) } ^{(1)}  :=  \sum_{v \in \mbox{\small pref} (\sigma(b))} |v|_b. 
\end{equation}
We assume that:
\begin{equation} \label{eqdelta}
\delta := \delta_{a  \in  \sigma(a) } ^{(1)}  - \delta_{a  \in  \sigma(b) } ^{(1)}  = 1.
\end{equation}

We fix for the rest of this work such a  substitution $\sigma$.
{}

\begin{lemma}
Let $\alpha\in\{a,b\}$ be a letter and $\chi_\alpha$ be the function defined for a word (finite or infinite) $v=v_0v_1\cdots$  by: 
\[
\chi_\alpha(v) =1  \mbox{ if }  v_0=\alpha \quad \mbox{and} \quad \chi_\alpha(v)=0   \mbox{ otherwise}.
\]
If $u=u_0u_1\cdots$ is a word (finite or infinite), we define \emph{the (left) shift} by $S(u)=u_1u_2\cdots$.
The terms $ \delta_{\alpha  \in  \sigma(a) } ^{(1)} $ are double Birkhoff sums:

\begin{equation} \label{eqdeltaabBis}
 \delta_{\alpha  \in  \sigma(a) } ^{(1)} 
 	=  \sum_{k=0}^{\lambda-1} (\lambda-k-1) \chi_\alpha \big(S^k \circ \sigma(a) \big)
 	=  \sum_{k=1}^{\lambda-1} \sum_{i=0} ^{k-1} \chi_\alpha \big(S^i \circ \sigma(a) \big)
	\quad \mbox{ for every $\alpha \in\{a,b\}$.}
\end{equation}
\end{lemma}
{}

\begin{definition} \label{defdeltaell}
We generalize Equation \eqref{eqdeltaabBis} to every positive integer $\ell$ by:
\begin{equation} \label{lesquatiazdi}
\delta_{\alpha  \in  \sigma(\beta)} ^{(\ell)} 
	:=\frac{1}{\ell!} \    \sum_{k=0}^{\lambda-1} (\lambda-k-1)^{\ell} \chi_\alpha \big(S^k \circ \sigma(\beta) \big)
\quad \mbox{for two letters $(\alpha,\beta) \in \{a,b\}^2$.}
\end{equation}
\end{definition}
{}

These terms are closed but different from the iterated Birkhoff sums over $\sigma(\alpha)$  introduced in \cite{BD} if $\ell > 2$.
By convention, we define $\delta_{a  \in  \sigma(a)} ^{(0)} :=\lambda_a$ and $\delta_{b  \in  \sigma(a)} ^{(0)} :=\lambda_b$.
It is clear that for every positive integer $\ell$ and every letter $\alpha\in\{a,b\}$:
\begin{equation} \label{azelkfqshglekz}
\delta_{a  \in  \sigma(\alpha)} ^{(\ell)}  + \delta_{b  \in  \sigma(\alpha)} ^{(\ell)}  =  \frac{1}{\ell!} \ \sum_{k=0}^{\lambda-1} (\lambda-k-1)^{\ell}.
\end{equation}
Note that it does not depend on the substitution.

\section{Definition of normalized moments}

\begin{definition} \label{defmoment}
Let $\sigma$ be a $\lambda$-uniform substitution satisfying  \eqref{eqlambda} and \eqref{eqdelta}.
Let $f$ be a solution of \eqref{eq:E_lambdq} which is the concatenation of two continuous functions $f_a,f_b:[0,\lambda] \rightarrow {\mathbb R}$ along a word $u = u_0u_1\cdots\in\{a,b\}^{\mathbb N}$ such that $u=\sigma(u)$.
We define the \emph{$\ell$-th moment} for $\ell \in {\mathbb N}$ by: 
\[
m_\alpha^{(\ell)} :=  \int_0^{\lambda} \Big( \lambda- x \Big) ^{\ell} \cdot f_\alpha(x)  d x \quad \mbox{for $\alpha \in \{a,b\}$.}
\]
We also define the \emph{$\ell$-th normalized moment} for $\ell \in {\mathbb N}$ by: 
\[
\tilde m_\alpha^{(\ell)} :=  \frac{1}{\ell!} \ \frac{1}{\lambda^ {\ell}} \   m_\alpha^{(\ell)}  = \frac{1}{\ell!} \ \frac{1}{\lambda^ {\ell}} \ \int_0^{\lambda} \Big( \lambda- x \Big) ^{\ell} \cdot f_\alpha(x)  d x
\quad \mbox{for $\alpha \in \{a,b\}$.}
\]
\end{definition}
{}

\begin{lemma} \label{leio1}
For every non-negative integer $\ell\in \mathbb N$ and every letter $\alpha \in \{a,b\}$:
\begin{equation} \label{leio2}
\sum_{q=0} ^{\ell+1}  (-1)^q 
\Big( \delta ^{(\ell+1-q)}_{a \in \sigma(\alpha)} \tilde m_a^{(q)} +  \delta ^{(\ell+1-q)}_{b \in \sigma(\alpha)} \tilde m_b^{(q)} \Big)
= - \frac{1}{\ell !} \lambda^{\ell+2} \ f \left( 0 \right) + \lambda^{\ell+1}	\tilde m_{\alpha}^{(\ell)}.
\end{equation}
\end{lemma}
{}

\noindent \emph{Proof of Lemma \ref{leio1}.}
We fix a letter $\alpha \in\{a,b\}$ and a non-negative integer $\ell \in \mathbb N$.
From Remark \ref{rem1}, the function $F(x) = \lambda \, f_{\alpha} \left( \frac x \lambda \right) $ is a primitive function of $f_{\sigma(\alpha)}$.
We calculate the following integral, recalling that $f_{\sigma(\alpha)}(0) = f_\alpha(0) = f(0)$  (Remark \ref{rem2}):
\begin{eqnarray}
 \int_0^{\lambda^2} \Big( \lambda^2- x \Big) ^{\ell+1} \cdot f_{\sigma(\alpha)}(x)  d x 
&=& \left[\Big( \lambda^2- x \Big) ^{\ell+1} \cdot F (x)   \right]_0^{\lambda^2} + (\ell+1) \int_0^{\lambda^2} \Big( \lambda^2- x \Big) ^{\ell} \cdot F (x)  d x \nonumber \\
&=&- \lambda^{2(\ell+1)} \cdot \lambda  \, f_{\alpha} \left( \frac 0 \lambda \right)+   (\ell+1) \int_0^{\lambda^2} \Big( \lambda^2- x \Big) ^{\ell} \cdot \lambda \,   f _{\alpha} \left( \frac x \lambda \right)  d x \nonumber \\
&=&	- \lambda^{2\ell+3} \cdot f(0) + (\ell+1) \cdot  \lambda^2 \int_0^{\lambda} \Big( \lambda^2 -  \lambda x \Big) ^{\ell}  \cdot  f_\alpha  \left( x  \right)  d x \nonumber \\
&=& - \lambda^{2\ell+3} \cdot f \left( 0 \right) + (\ell+1) \cdot	\lambda^{\ell+2} m_\alpha^{(\ell)} \label{zqfgzhek}.
\end{eqnarray}

We write $\sigma(\alpha)=v_0\cdots v_{\lambda-1}$ and with Definition \ref{conca}:
\begin{align*}
\int_0^{\lambda^2} \Big( \lambda^2- x \Big) ^{\ell+1} \cdot f_{\sigma(\alpha)}(x)  d x
&= \sum_{k=0}^{\lambda-1} \int_{k\lambda} ^{(k+1)\lambda} \Big( \lambda^2- x \Big) ^{\ell+1} \cdot f_{\sigma(\alpha)}(x)  d x \\
&= \sum_{k=0}^{\lambda-1} \int_{0} ^{\lambda} \Big( \lambda^2- x -k\lambda \Big) ^{\ell+1} \cdot f_{v_k} (x)  d x.
\end{align*}

It  remains to express $( \lambda^2- x -k\lambda) ^{\ell+1} $  in the basis $((\lambda-x)^q;0\leq q \leq \ell+1)$.
To do this, we derive $q$ times the polynomial function $( \lambda^2- x -k\lambda ) ^{\ell+1} $ and we estimate it at $x = \lambda$.
We fix $k\in\{0,\ldots,\lambda-1\}$:
\[
\Big( \lambda^2- x -k\lambda  \Big) ^{\ell+1}   = \sum_{q=0} ^{\ell+1} a_{q,k} ^{(\ell+1)} (\lambda-x)^{q} ,
\]
where $a_{q,k}^{(\ell+1)}   = \frac{(\ell+1)!}{q! \ (\ell+1-q)!} \ (-1)^{q} \ \Big(\lambda^{2}-\lambda(k+1)  \Big)^{\ell+1-q}$
for $0\leq q \leq \ell+1$.
In particular:
\[
a_{\ell+1,k}^{(\ell+1)}  =   (-1)^{\ell+1} , \  a_{\ell,k}^{(\ell+1)}  =  (\ell+1) \cdot (-1)^{\ell }\cdot \lambda \cdot \Big(\lambda-k-1  \Big), \ \ldots
\]

Therefore we have:
\begin{equation} \label{zqfgzhekBis}
\int_0^{\lambda^2}	\Big( \lambda^2- x \Big) ^{\ell+1} \cdot f_{\sigma(\alpha)}(x)  d x =
\sum_{k=0}^{\lambda-1} \sum_{q=0} ^{\ell+1}  a_{q,k}^{(\ell+1)} 	  \int_{0} ^{\lambda}\Big( \lambda- x \Big) ^q \cdot f_{v_k} (x)  d x.
\end{equation}

With Equations \eqref{zqfgzhek}, we find:
\[
\sum_{q=0} ^{\ell+1}  \sum_{k=0}^{\lambda-1} a_{q,k}^{(\ell+1)}    m_{v_k}^{(q)}  = - \lambda^{2\ell+3} \ f \left( 0 \right) + (\ell+1)	\lambda^{\ell+2} m_{\alpha}^{(\ell)}.
\]

The normalized relation is:
\[
\sum_{q=0} ^{\ell+1}  \frac{(-1)^q }{(\ell+1-q)!}  \sum_{k=0}^{\lambda-1}  \Big(\lambda-k-1  \Big)^{\ell+1-q}     \tilde m_{v_k}^{(q)} 
= - \frac{1}{\ell !} \lambda^{\ell+2} \ f \left( 0 \right) + \frac 1 {\ell !} 	\lambda	m_{\alpha}^{(\ell)}.
\]

We simplify this expression with Definition \ref{defdeltaell} of $ \delta ^{(\ell+1-q)}_{\beta \in \sigma(\alpha)}$: \smallskip

~\hfill
$\displaystyle \sum_{q=0} ^{\ell+1}  (-1)^q  \Big( \delta ^{(\ell+1-q)}_{a \in \sigma(\alpha)} \tilde m_a^{(q)} +  \delta ^{(\ell+1-q)}_{b \in \sigma(\alpha)} \tilde m_b^{(q)} \Big)
= - \frac{1}{\ell !} \lambda^{\ell+2} \ f \left( 0 \right) + \lambda^{\ell+1}	\tilde m_{\alpha}^{(\ell)}.$ ~
\hfill ~\qed

\section{A technical lemma}

\begin{lemma} \label{lemma:uniciteA}
Let $\lambda$ be a positive real number and $f$ be a continuous function which is solution of \eqref{eq:E_lambdq}. 
Then for every  $n\in {\mathbb N}$ and $\ell\in\mathbb N$:
\begin{equation} \label{equationuniciteA}
f\left(\frac{n+1}{\lambda^{\ell+1}}\right) -f\left(\frac{n}{\lambda^{\ell+1}}\right)   	= 
	\sum \limits_{k=0} ^{\ell} \frac{1}{(k+1)!}  \, \frac{1}{\lambda^{(k+1)\left( \ell- \frac{k}{2} \right)} } f\left(\frac{n}{\lambda^{\ell-k}}\right)  + I_n^{(\ell+1)},
\end{equation}
where the remainder integral  $I_n^{(\ell)}$ is:
\begin{equation} \label{equadefresteintegrlae}
I_n^{(\ell)} :=  \frac{1}{\lambda^{\ell(\ell+1)/2}} \, \frac{1}{\ell!} \int_0^\lambda   (\lambda-u)^{\ell} f(u+\lambda n)   d u.
\end{equation}
\end{lemma}
{}

\begin{remark}
Let $\sigma$ be a $\lambda$-uniform substitution satisfying  \eqref{eqlambda} and \eqref{eqdelta}.
We suppose that $f$ is a solution of \eqref{eq:E_lambdq} which is the concatenation of two continuous functions $f_a,f_b:[0,\lambda] \rightarrow {\mathbb R}$ along a word $u = u_0u_1\cdots\in\{a,b\}^{\mathbb N}$ such that $u=\sigma(u)$.

Then for every $n \in {\mathbb N}$ and $\ell\in {\mathbb N}$,  $I_n^{(\ell)}$  depends only on $u_n$ and $\ell$. 
With Definition \ref{defmoment} of moments, Relation \eqref{equationuniciteA} can be rewritten as follows:
\begin{equation} \label{equationuniciteABis}
f\left(\frac{n+1}{\lambda^{\ell+1}}\right) -f\left(\frac{n}{\lambda^{\ell+1}}\right)   	= 
	\sum \limits_{k=0} ^{\ell} \frac{1}{(k+1)!}  \, \frac{1}{\lambda^{(k+1)\left( \ell- \frac{k}{2} \right)} }\  f\left(\frac{n}{\lambda^{\ell-k}}\right)  +  \frac{1}{\lambda^{\ell(\ell-1)/2}} \ \tilde m_{v_n}^{(\ell)}.
\end{equation}
\end{remark}
{}

\noindent \emph{Proof of Lemma \ref{lemma:uniciteA}.}
We fix two non-negative integers $n$ and $\ell$. From Equation \eqref{eq:E_lambdq}:
\[ f\left(\frac{n+1}{\lambda^{\ell+1}}\right)  -    f\left(\frac{n}{\lambda^{\ell+1}}\right) = \int_{n/\lambda^\ell}^{(n+1)/\lambda^\ell} f(t) d t. \]

Still according to Equation \eqref{eq:E_lambdq}, 
the values of the function at $n\cdot  \lambda^{-\ell}$ and $t\in {\mathbb R}_+$ satisfy:
\[ \forall t \in {\mathbb R}_+, \quad f(t) =  f\left(\frac{n}{\lambda^\ell}\right) + \int_{n/\lambda^{\ell-1}} ^{\lambda t} f(s_1)  d s_1. \]

The two previous relations involve:
\[
f\left(\frac{n+1}{\lambda^{\ell+1}}\right)  -    f\left(\frac{n}{\lambda^{\ell+1}}\right) 
= \frac{ 1}{\lambda^{\ell}}f\left(\frac{n}{\lambda^\ell}\right) + \int_{n/\lambda^\ell}^{(n+1)/\lambda^\ell}   \int_{n/\lambda^{\ell-1}} ^{\lambda t} f(s_1)  d s_1   d t .
\]

We can iterate the process using the relation:
\[ \mbox{for $s_1 \in {\mathbb R}_+$,} \quad f(s_1) = f\left(\frac{n}{\lambda^{\ell-1}}\right) +\int_{n/\lambda^{\ell-2}} ^{\lambda s_1} f(s_2)  d s_2 . \]

The goal is to continue this process (like Taylor series) 
and to express $f\left(\frac{n+1}{\lambda^{\ell+1}}\right)  -    f\left(\frac{n}{\lambda^{\ell+1}}\right)$
as a linear combinaison of $\{f(n\cdot \lambda^{-k}); 0 \leq k \leq \ell \}$
and an integral on $[k\lambda,(k+1)\lambda]$.
To do this, we introduce the following terms:
\begin{align}
&V_k :=\int_{n/\lambda^\ell}^{(n+1)/\lambda^\ell}	\int_{n/\lambda^{\ell-1}}^{\lambda t}	\int_{n/\lambda^{\ell-2}}^{\lambda s_1}	\cdots	\int_{n/\lambda^{\ell-k}}^{\lambda s_{k-1}}   d s_k\cdots  d s_1  d t
\quad \mbox{for $k\in\{0,\ldots,\ell\}$,}   \label{eqkjfblqsk} \\
&I := 	\int_{n/\lambda^\ell}^{(n+1)/\lambda^\ell}   \int_{n/\lambda^{\ell-1}}^{\lambda t} \int_{n/\lambda^{\ell-2}}^{\lambda s_1} \cdots  \int_{\lambda n}^{\lambda s_{\ell}}  f(s_{\ell+1})  d s_{\ell+1} \cdots  d s_1   d t .
\label{eqkjfblqskbis} 
\end{align}

Continuing the process described above, we have:
\begin{equation} \label{bnqcvcqnvnq1}
f\left(\frac{n+1}{\lambda^{\ell+1}}\right)  -    f\left(\frac{n}{\lambda^{\ell+1}}\right)  = \sum_{k=0}^{\ell} V_{k}  f \left(\frac{n}{\lambda^{\ell-k}}\right) + I.
\end{equation} 
It remains to calculate the values of $(V_k)_{0\leq k \leq \ell}$ and $I$. \medskip

\begin{itemize}
\item 
We fix $k\in\{0,\ldots,\ell\}$ and we calculate $V_k$ defined in \eqref{eqkjfblqsk}.
We put  $u_i = \lambda (\lambda^{\ell-i} s_{i}-n )  $ for $i\in \{1,\ldots,k-1\}$ and $v= \lambda(\lambda ^\ell t-n) $:
\begin{eqnarray}
V_k   
&=&  \prod \limits_{i=0}^{k} \frac{1}{\lambda^{\ell-i+1}} \int_0^\lambda  \int_0^v \int_0^{u_1} \cdots \int_0^{u_{k-1}}  d u_k \cdots  d u_1 d v. 
	\nonumber \\
&=& \frac{1}{\lambda^{(k+1)(\ell+1)}}  \left( \prod \limits_{i=0}^{k}  \lambda^{i}\right) \frac{\lambda^{k+1}}{(k+1)!} =  \frac{1}{(k+1)!}   \lambda^{(k+1)\left( \frac{k}{2}-\ell \right)}.  \label{bnqcvcqnvnq2}
\end{eqnarray}

\item 
We make the substitution $u+\lambda n = s_{\ell+1}$ in Equation \eqref{eqkjfblqskbis}:
\[
I = 	\int_{n/\lambda^\ell}^{(n+1)/\lambda^\ell}   \int_{k/\lambda^{\ell-1}}^{\lambda t} \int_{n/\lambda^{\ell-2}}^{\lambda s_1} 
\cdots  \int_{\lambda n}^{\lambda s_{\ell}}  f (u+\lambda n )  d u  d s_{\ell} \cdots  d s_1   d t.
\]

If we put $u_i = \lambda (\lambda^{\ell-i} s_{i}-n )$ for $i\in\{1,\ldots,\ell\}$ and $v= \lambda(\lambda ^\ell t-n)$:
\begin{eqnarray}
I 
&=&  \prod \limits_{i=0}^{\ell+1} \frac{1}{\lambda^{\ell-i+1}} \int_0^\lambda  \int_0^v \int_0^{u_1}\cdots \int_0^{u_{\ell-1}}  f (u+\lambda n )   d u  d u_\ell \cdots  d u_1 d v \nonumber \\
&= &\prod \limits_{i=0}^{\ell+1} \frac{1}{\lambda^{i}} \int_0^\lambda   \left( \int_{u<u_{\ell}<\cdots<u_1<v<\lambda}  d v d u_1\cdots  d u_\ell \right)  f (u+\lambda n )   d u \nonumber \\
&= &\frac{1}{\lambda^{(\ell+1)(\ell+2)/2}} \int_0^\lambda   \frac{(\lambda-u)^{\ell+1}}{(\ell+1)!} f (u+\lambda n )   d u = I_n^{(\ell+1)}. \label{bnqcvcqnvnq3} 
\end{eqnarray}
\end{itemize}
Lemma \ref{lemma:uniciteA} is proved by combining Relations \eqref{bnqcvcqnvnq1}, \eqref{bnqcvcqnvnq2} and \eqref{bnqcvcqnvnq3}. ~ \hfill ~ \qed

\section{Proof of Theorem \ref{theorem}} \label{proof}

Let $\sigma$ be a $\lambda$-uniform substitution satisfying  \eqref{eqlambda} and \eqref{eqdelta}.
We denote by $\mathcal S_\lambda$ the set of solutions of \eqref{eq:E_lambdq} which are the concatenation of two continuous functions $[0,\lambda] \rightarrow {\mathbb R}$ along a word $u = u_0u_1\cdots\in\{a,b\}^{\mathbb N}$ such that $u=\sigma(u)$. \medskip

\emph{We prove that the map from $\mathcal S_\lambda$ into $\mathbb R^2$ defined by $f\mapsto \big(f(0),f(1)\big)$ is an injective morphism.} \medskip

We fix a function $f\in \mathcal S_\lambda$  such that $f(0)=f(1)=0$.

\begin{enumerate}
\item 
From Equation \eqref{eq:E_lambdq}:
\[ m_{u_0}^{(0)} = \tilde m_{u_0}^{(0)} = \int_0^\lambda  f(t)  d t = f(1)-f(0)=0.  \]
We calculate the following integral for every non-negative integer $n$:
\[
\int_0^{n \lambda^2}  f(t)  d t = 
\left\{ \begin{array}{l}
\displaystyle  \sum_{k=0}^{n-1}\sum_{i=0}^{\lambda-1}  \int_{k^2\lambda+i\lambda} ^{k^2\lambda+(i+1)\lambda}  f(t)  d t 
	= \sum_{k=0}^{n-1} \Big(\lambda_a m_a^{(0)} + \lambda_b m_b^{(0)} \Big)=n \Big(\lambda_a m_a^{(0)} + \lambda_b m_b^{(0)} \Big), \\
f(n\lambda)-f(0).
\end{array} \right.
\]
We  divide this expression by $n$ and since $f$ is bounded:
\[  \lambda_a m_a^{(0)} + \lambda_b m_b^{(0)} = \frac 1 nf(n\lambda)-\frac 1 n f(0) \underset{n \to +\infty}{\longrightarrow} 0 . \]
So we have $m_a^{(0)}=m_b^{(0)}=0$ and for every $n \in \mathbb N$:
\begin{align*}
f(n) 
	&= f(0)+\int_0^{n\lambda} f(t)  d t = f(0) + \sum_{k=0}^{n-1} \int_0^\lambda f(t)  d t
	= f(0) + \sum_{k=0}^{n-1} m_{u_k}^{(0)} \\
	&= f(0) + |u_0\cdots u_{n-1}|_a \cdot m_{a}^{(0)}+ |u_0\cdots u_{n-1}|_b  \cdot m_{b}^{(0)} = 0.
\end{align*}

\item \emph{We show by induction on $\ell\geq 1$ that $m_{\alpha}^{(i)}=0$ for every $\alpha \in \{a,b\}$ and every $0\leq i \leq \ell$.} \medskip 

We suppose that $m_{\alpha}^{(i)}=0$ for $\alpha \in \{a,b\}$ and $0\leq i \leq \ell$.
We recall Relation \eqref{leio2} in Lemma \ref{leio1} for $\alpha=a$:
\[
\sum_{q=0} ^{\ell+1}  (-1)^q  \Big( \delta ^{(\ell+1-q)}_{a \in \sigma(a)} \tilde m_a^{(q)} +  \delta ^{(\ell+1-q)}_{b \in \sigma(a)} \tilde m_b^{(q)} \Big)
= - \frac{1}{\ell !} \lambda^{\ell+2} \ f \left( 0 \right) + \lambda^{\ell+1}	\tilde m_{a}^{(\ell)}.
\]

By induction hypothesis, we have:
\begin{equation} \label{premiererelat}
\lambda_a \tilde m_a^{(\ell+1)} + \lambda_b \tilde m_b^{(\ell+1)} = 0.
\end{equation}

We use Relation \eqref{leio2}  for the positive integer $\ell+1$ with $\alpha=a$ and $\alpha = b$:
\[
\begin{tabular}{cl}
&
$ \left\{ \begin{array}{l}
\displaystyle \sum_{q=0} ^{\ell+2}  (-1)^q  \Big( \delta ^{(\ell+2-q)}_{a \in \sigma(a)} \tilde m_a^{(q)} +  \delta ^{(\ell+2-q)}_{b \in \sigma(a)} \tilde m_b^{(q)} \Big)
= - \frac{1}{(\ell+1) !} \lambda^{\ell+3} \ f \left( 0 \right) + \lambda^{\ell+2}	\tilde m_{a}^{(\ell+1)}, \\
\displaystyle \sum_{q=0} ^{\ell+2}  (-1)^q \Big( \delta ^{(\ell+2-q)}_{a \in \sigma(b)} \tilde m_a^{(q)}  +  \delta ^{(\ell+2-q)}_{b \in \sigma(b)} \tilde m_b^{(q)} \Big)
= - \frac{1}{(\ell+1) !} \lambda^{\ell+3} \ f \left( 0 \right) + \lambda^{\ell+2}	\tilde m_{b}^{(\ell+1)}.
\end{array} \right.$
\\
$\iff$ &
$ \left\{ \begin{array}{l}
\displaystyle \sum_{q=\ell+1} ^{\ell+2}  (-1)^q  \Big( \delta ^{(\ell+2-q)}_{a \in \sigma(a)} \tilde m_a^{(q)} +  \delta ^{(\ell+2-q)}_{b \in \sigma(a)} \tilde m_b^{(q)} \Big) =  \lambda^{\ell+2}	\tilde m_{a}^{(\ell+1)}, \\
\displaystyle \sum_{q=\ell+1} ^{\ell+2}  (-1)^q \Big( \delta ^{(\ell+2-q)}_{a \in \sigma(b)} \tilde m_a^{(q)}  +  \delta ^{(\ell+2-q)}_{b \in \sigma(b)} \tilde m_b^{(q)} \Big)=  \lambda^{\ell+2}	\tilde m_{b}^{(\ell+1)}.
\end{array} \right. $
\end{tabular}
\]

When we subtract these two relations, the coefficients of $\tilde m_a^{(\ell+2)}$ and $\tilde m_a^{(\ell+2)}$ vanish:
\begin{equation} \label{sagrdaztta}
(-1)^{\ell+1} \cdot \Big( \delta ^{(1)}_{a \in \sigma(a)} -  \delta ^{(1)}_{a \in \sigma(b)} \Big) \tilde m_a^{(\ell+1)}  
	+ (-1)^{\ell+1} \cdot \Big( \delta ^{(1)}_{b \in \sigma(a)} -  \delta ^{(1)}_{b \in \sigma(b)} \Big) \tilde m_a^{(\ell+1)}  =  \lambda^{\ell+2}	\Big( \tilde m_{a}^{(\ell+1)} - \tilde m_{b}^{(\ell+1)} \Big).
\end{equation}
We recall that $\delta ^{(1)}_{a \in \sigma(a)} -  \delta ^{(1)}_{a \in \sigma(b)}=1$. Moreover from \eqref{azelkfqshglekz},
\[ \delta ^{(1)}_{a \in \sigma(a)}+\delta ^{(1)}_{b \in \sigma(a)} =  \delta ^{(1)}_{a \in \sigma(b)}  + \delta ^{(1)}_{b \in \sigma(b)} . \]
Then  $\delta ^{(1)}_{a \in \sigma(b)} -  \delta ^{(1)}_{a \in \sigma(a)}=-1$. We can rewrite \eqref{sagrdaztta}:
\[
(-1)^{\ell+1}  \tilde m_a^{(\ell+1)}  - (-1)^{\ell+1}  \tilde m_a^{(\ell+1)} =  \lambda^{\ell+2}	\Big( \tilde m_{a}^{(\ell+1)} - \tilde m_{b}^{(\ell+1)} \Big),
\]
so we find:
\begin{equation} \label{premiererelatBis}  \tilde m_a^{(\ell+1)} - \tilde m_b^{(\ell+1)} = 0. \end{equation}

Combining Equations \eqref{premiererelat} and \eqref{premiererelatBis}: $\tilde m_a^{(\ell+1)}=\tilde m_b^{(\ell+1)}=0$. 
\bigskip

\item \emph{We show by induction on $\ell\in {\mathbb N}$ that $f(n \lambda^{-i})=0$  for every $n \in {\mathbb N}$ and $0\leq i \leq \ell$.} \medskip 

For every $n\in \mathbb N$, we recall Equation \eqref{equationuniciteABis}:
\[
f\left(\frac{n+1}{\lambda^{\ell+1}}\right) -f\left(\frac{n}{\lambda^{\ell+1}}\right)   	= \sum \limits_{k=0} ^{\ell} \frac{1}{(k+1)!}  \, \frac{1}{\lambda^{(k+1)\left( \ell- \frac{k}{2} \right)} }\  f\left(\frac{n}{\lambda^{\ell-k}}\right).
\]
It is easy to verify by induction that for every $n\in \mathbb N$:
\[
f\left(\frac{n+1}{\lambda^{\ell+1}}\right) -f\left(\frac{n}{\lambda^{\ell+1}}\right)   	= 0.
\]
Since $f(0)=0$, therefore $f\left(\frac{n}{\lambda^{\ell+1}}\right)   =0$ for every $n\in \mathbb N$.

\item
We have seen that $f$  vanishes at $\lambda$-adic points  (i.e. the points $n \lambda^{-k}$ for $n,k\in {\mathbb N}$).
These points form a dense subset of ${\mathbb R}_+$ and $f$ is a continuous function, therefore $f$ is the identically zero function.
\end{enumerate}




\end{document}